\documentclass[11pt]{article}

\usepackage{amssymb, amsmath, amsthm, graphicx}
\usepackage[left=1in,top=1in,right=1in]{geometry}

\usepackage{subfigure}

      \theoremstyle{plain}
      \newtheorem{theorem}{Theorem}[section]
      \newtheorem{lemma}[theorem]{Lemma}
            
            \newtheorem{observation}[theorem]{Observation}

      \theoremstyle{definition}

      \theoremstyle{remark}

\title{The number of edges in $k$-quasi-planar graphs\footnote{A preliminary version of this paper with A. Suk as its sole author will appear in {\em Proc. 19th Internat. Symp. on Graph Drawing (GD 2011, TU Eindhoven), LNCS}, Springer, 2011}}

\author{Jacob Fox \thanks{Massachusetts Institute of Technology, Cambridge, MIT. Supported by a Simons Fellowship and NSF grant DMS 1069197. Email: {\tt fox@math.mit.edu}} \and J\'anos Pach\thanks{EPFL, Lausanne and Courant Institute,  New York, NY. Supported by Hungarian Science Foundation EuroGIGA Grant OTKA NN 102029, by Swiss National Science Foundation Grant 200021-125287/1, and by NSF Grant CCF-08-30272. Email: {\tt pach@cims.nyu.edu}.} \and Andrew Suk\thanks{Massachusetts Institute of Technology, Cambridge. Supported by an NSF Postdoctoral Fellowship. Email: {\tt asuk@math.mit.edu}.}
}

\begin{document}

\maketitle

\begin{abstract}

A graph drawn in the plane is called \emph{$k$-quasi-planar} if it does not contain $k$ pairwise crossing edges. It has been conjectured for a long time that for every fixed $k$, the maximum number of edges of a $k$-quasi-planar graph with $n$ vertices is $O(n)$. The best known upper bound is $n(\log n)^{O(\log k)}$. In the present note, we improve this bound to $(n\log n )2^{\alpha^{c_k}(n)}$ in the special case where the graph is drawn in such a way that every pair of edges meet at most once. Here $\alpha(n)$ denotes the (extremely slowly growing) inverse of the Ackermann function. We also make further progress on the conjecture for $k$-quasi-planar graphs in which every edge is drawn as an $x$-monotone curve. Extending some ideas of Valtr, we prove that the maximum number of edges of such graphs is at most $2^{ck^6}n\log n$.

\end{abstract}

\section{Introduction}

A \emph{topological graph} is a graph drawn in the plane such that its vertices are represented by points and its edges are represented by non-self-intersecting arcs connecting the corresponding points. In notation and terminology, we make no distinction between the vertices and edges of a graph and the points and arcs representing them, respectively. No edge is allowed to pass through any point representing a vertex other than its endpoints. Any two edges can intersect only in a finite number of points. Tangencies between the edges are not allowed. That is, if two edges share an interior point, then they must properly cross at this point. A topological graph is \emph{simple} if every pair of its edges intersect at most once: at a common endpoint or at a proper crossing. If the edges of a graph are drawn as straight-line segments, then the graph is called \emph{geometric}.

Finding the maximum number of edges in a topological graph with a forbidden crossing pattern is a fundamental problem in extremal topological graph theory (see \cite{ackerman2,acktardos,agarwal,cap,fox,fulek,pinchasi,tardos,valtr}).  It follows from Euler's Polyhedral Formula that every topological graph on $n$ vertices and with no two crossing edges has at most $3n-6$ edges.  A graph is called $k$\emph{-quasi-planar} if it can be drawn as a topological graph with no $k$ pairwise crossing edges.  A graph is 2-quasi-planar if and only if it is planar.  According to an old conjecture (see Problem 1 in Section 9.6 of \cite{brass}), for any fixed $k \ge 2$ there exists a constant $c_k$ such that every $k$-quasi-planar graph on $n$ vertices has at most $c_kn$ edges.  Agarwal, Aronov, Pach, Pollack, and Sharir \cite{agarwal} were the first to prove this conjecture for {\em simple} 3-quasi-planar graphs.  Later, Pach, Radoi\v{c}i\'c, and T\'oth \cite{rados} generalized the result to {\em all} 3-quasi-planar graphs.  Ackerman \cite{ackerman} proved the conjecture for $k = 4$.

For larger values of $k$, first Pach, Shahrokhi, and Szegedy \cite{pach} showed that every simple $k$-quasi-planar graph on $n$ vertices has at most $c_kn(\log n)^{2k-4}$ edges. For $k\ge 3$ and for all (not necessarily simple) $k$-quasi-planar graphs, Pach, Radoi\v{c}i\'c, and T\'oth \cite{rados} established the upper bound $c_kn(\log n)^{4k-12}$. Plugging into these proofs the above mentioned result of Ackerman \cite{ackerman}, for $k\ge 4$, we obtain the slightly better bounds $c_kn(\log n)^{2k-8}$ and $c_kn(\log n)^{4k-16}$, respectively. For large values of $k$, the exponent of the polylogarithmic factor in these bounds was improved by Fox and Pach \cite{fox}, who showed that the maximum number of edges of a $k$-quasi-planar graph on $n$ vertices is $n(\log n)^{O(\log k)}$.

For the number of edges of geometric graphs, that is, graphs drawn by straight-line edges, Valtr \cite{valtrpar} proved the upper bound $O(n\log n)$.  He also extended this result to {\em simple} topological graphs whose edges are drawn as {\em $x$-monotone} curves \cite{valtr}.

The aim of this paper is to improve the best known bound, $n(\log n)^{O(\log k)}$, on the number of edges of a $k$-quasi-planar graph in two special cases: for simple topological graphs and for not necessarily simple topological graphs whose edges are drawn as $x$-monotone curves. In both cases, we improve the exponent of the polylogarithmic factor from $O(\log k)$ to $1+o(1)$.

\begin{theorem}
\label{simple}
Let $G = (V,E)$ be a $k$-quasi-planar simple topological graph with $n$ vertices.  Then $|E(G)| \leq   (n\log n)2^{\alpha(n)^{c_k}}$, where $\alpha(n)$ denotes the inverse of the Ackermann function and $c_k$ is a constant that depends only on $k$.
\end{theorem}

Recall that the {\em Ackermann} (more precisely, the {\em Ackermann-P\'eter) function} $A(n)$ is defined as follows.  Let $A_1(n) = 2n$, and $A_k(n) = A_{k-1}(A_k(n-1))$ for $k=2,3,\ldots$. In particular, we have $A_2(n) = 2^n$, and $A_3(n)$ is an exponential tower of $n$ {\em two}'s. Now let $A(n) = A_n(n)$, and let $\alpha(n)$ be defined as $\alpha(n) = \min\{k \geq 1: A(k) \geq n\}$. This function grows much slower than the inverse of any primitive recursive function.

\begin{theorem}
\label{xmono}
Let $G =(V,E)$ be a $k$-quasi-planar (not necessarily simple) topological graph with $n$ vertices, whose edges are drawn as $x$-monotone curves.  Then $|E(G)| \leq 2^{ck^{6}}n\log n$, where $c$ is an absolute constant.
\end{theorem}

In both proofs, we follow the approach of Valtr \cite{valtr} and apply results on generalized Davenport-Schinzel sequences. An important new ingredient of the proof of Theorem \ref{simple} is a corollary of a separator theorem established in \cite{fox3} and developed in \cite{FoPa}. Theorem \ref{xmono} is not only more general than Valtr's result, which applies only to simple topological graphs, but its proof gives a somewhat better upper bound: we use a result of Pettie \cite{seth}, which improves the dependence on $k$ from double exponential to single exponential.

\section{Generalized Davenport-Schinzel Sequences}

The sequence $u = a_1,a_2,...,a_m$ is called {\em $l$-regular} if any $l$ consecutive terms are pairwise different.  For integers $l,t \geq 2$, the sequence

$$S = s_1,s_2,...,s_{l t}$$

\noindent of length $lt$ is said to be of {\em type $up(l,t)$} if the first $l$ terms are pairwise different and

$$s_i = s_{i +l } = s_{i + 2l} = \cdots = s_{i + (t-1)l }$$

\noindent for every $i, 1\le i\le l$. For example,

$$a,b,c,a,b,c,a,b,c,a,b,c,$$

\noindent is a type $up(3,4)$ sequence or, in short, an $up(3,4)$ sequence. We need the following theorem of Klazar \cite{klazar92} on generalized Davenport-Schinzel sequences.

\begin{theorem}[Klazar]
\label{klazar}
For $l\geq 2$ and $t \geq 3$, the length of any $l$-regular sequence over an $n$-element alphabet that does not contain a subsequence of type $up(l,t)$ has length at most

$$n\cdot l 2^{(lt-3)} \cdot (10l)^{10\alpha(n)^{lt}}.$$

\end{theorem}

\noindent  For $l \geq 2$, the sequence

$$S = s_1,s_2,...,s_{3l-2}$$

\noindent of length $3l-2$ is said to be of {\em type up-down-up}$(l)$ if the first $l$ terms are pairwise different and
$$s_i = s_{2l-i} = s_{(2l-2)+ i}$$

\noindent for every $i, 1\le i\le i$. For example,

$$a,b,c,d,c,b,a,b,c,d,$$

\noindent is an \emph{up-down-up}$(4)$ sequence.  Valtr and Klazar \cite{val} showed that any $l$-regular sequence over an $n$-element alphabet, which contains no subsequence of type up-down-up$(l)$, has length at most $2^{l^c}n$ for some constant $c$.  This has been  improved by Pettie \cite{seth}, who proved the following.

\begin{lemma}[Pettie]
\label{updownup}
For $l\geq 2$, the length of any $l$-regular sequence over an $n$-element alphabet, which contains no subsequence of type up-down-up$(l)$, has length at most $2^{O(l^2)}n$.
\end{lemma}

\noindent For more results on generalized Davenport-Schinzel sequences, see \cite{gabriel,seth,pettie}.

\section{On intersection graphs of curves}
In this section, we prove a useful lemma on intersection graphs of curves. It shows that every collection $C$ of curves, no two of which intersect many times,  contains a large subcollection $C'$ such that in the partition of $C'$ into its connected components $C_1,\ldots,C_t$ in the intersection graph of $C$, each component $C_i$ has a vertex connected to all other $|C_i|-1$ vertices.

For a graph $G=(V,E)$, a subset $V_0$ of the vertex set is said to be a {\em separator} if there is a partition $V=V_0 \cup V_1 \cup V_2$ with $|V_1|,|V_2| \leq \frac{2}{3}|V|$ such that no edge connects a vertex in $V_1$ to a vertex in $V_2$. We need the following separator lemma for intersection graphs of curves, established in \cite{fox3}.

\begin{lemma}[Fox--Pach]
 \label{seplem}
There is an absolute constant $c_1$ such that
every collection $C$ of curves with $x$ intersection points has a separator of size at most $c_1\sqrt{x}$.
\end{lemma}

Call a collection $C$ of curves in the plane {\it decomposable} if there is a partition $C=C_1 \cup \ldots \cup C_t$ such that each $C_i$ contains a curve which intersects all other curves in $C_i$, and for $i \not = j$, the curves in $C_i$ are disjoint from the curves in $C_j$. The following lemma is a strengthening of Proposition 6.3 in \cite{FoPa}. Its proof is essentially the same as that of the original statement. It is include here, for completeness.

\begin{lemma}
\label{decompose}
There is an absolute constant $c>0$ such that every collection $C$ of $m \geq 2$ curves such that each pair of them intersect in at most $t$ points has a decomposable subcollection of size at least  $\frac{cm}{t\log m}$.
\end{lemma}

\noindent{\bf Proof of Lemma \ref{decompose}} We prove the following stronger statement. There is an absolute constant $c>0$ such that every collection $C$ of $m \geq 2$ curves whose intersection graph has at least $x$ edges, and each pair of curves intersect in at most $t$ points, has a decomposable subcollection of size at least  $\frac{cm}{t\log m}+\frac{x}{m}$.  Let $c=\frac{1}{576c_1^2}$, where $c_1 \geq 1$ is the constant in Lemma \ref{seplem}. The proof is by induction on $m$, noting that all collections of curves with at most three elements are decomposable. Define $d=d(m,x,t):=\frac{cm}{t\log m}+\frac{x}{m}$.

Let $\Delta$ denote the maximum degree of the intersection graph of $C$. We have $\Delta < d-1$. Otherwise, the subcollection consisting of a curve of maximum degree, together with the curves in $C$ that intersect it, is decomposable and its size is at least $d$, and we are done. Also, $\Delta \geq 2\frac{x}{m}$, since $2\frac{x}{m}$ is the average degree of the vertices in the intersection graph of $C$. Hence, if $\Delta \geq 2\frac{cm}{t\log m}$, then the desired inequality holds. Thus, we may assume $\Delta < 2\frac{cm}{t\log m}$.

Applying Lemma \ref{seplem} to the intersection graph of $C$, we obtain that there is a separator $V_0 \subset C$ with $|V_0| \leq c_1\sqrt{tx}$, where $c_1$ is the absolute constant in Lemma \ref{seplem}. That is, there is a partition $C=V_0 \cup V_1 \cup V_2$ with $|V_1|,|V_2| \leq 2|V|/3$ such that no curve in $V_1$ intersects any curve in $V_2$. For $i=1,2$, let $m_i=|V_i|$ and $x_i$ denote the number of pairs of curves in $V_i$ that intersect, so that
\begin{equation}\label{equ1}
x_1+x_2 \geq x-\Delta|V_0| \geq x-2\frac{cm}{t\log m}c_1\sqrt{tx}.\end{equation}
As no curve in $V_1$ intersects any curve in $V_2$, the union of a decomposable subcollection of $V_1$ and a decomposable subcollection of $V_2$ is decomposable. Thus, by the induction hypothesis, $C$ contains decomposable subcollection of size at least
\begin{eqnarray*} d(m_1,x_1,t)+d(m_2,x_2,t) & = & \frac{cm_1}{t\log m_1}+\frac{x_1}{m_1}+\frac{cm_2}{t\log m_2}+\frac{x_2}{m_2} \\ & \geq &  \frac{c(m_1+m_2)}{t\log (2m/3)}+\frac{(x_1+x_2)}{2m/3}
.\end{eqnarray*}
We split the rest of the proof into two cases.

\noindent \emph{Case 1.} $x \geq t^{-1}\left(12c_1c\frac{m}{\log m}\right)^2$. In this case, by (\ref{equ1}), we have $x_1+x_2 \geq \frac{5}{6}x$ and hence there is a decomposable subcollection of size at least
\begin{eqnarray*} d(m_1,x_1,t) + d(m_2,x_2,t)  & \geq & \frac{c(m_1+m_2)}{t\log m}+\frac{5x}{4m}
  =  d+\frac{x}{4m}-\frac{c(m-(m_1+m_2))}{t\log m} \\ & \geq &
d+\frac{x}{4m}-\frac{c_1 c\sqrt{tx}}{t\log m} > d,
\end{eqnarray*}
completing the analysis.

\noindent \emph{Case 2.} $x < t^{-1}\left(12c_1c\frac{m}{\log m}\right)^2$. There is a decomposable subcollection of size at least
\begin{eqnarray*}
d(m_1,x_1,t) + d(m_2,x_2,t) & \geq & \frac{c(m_1+m_2)}{t\log (2m/3)}  \geq  \frac{c}{t}\left(m-c_1\sqrt{tx}\right)\left(\frac{1}{\log m}+\frac{1}{2\log^2 m}\right) \\ & \geq & \frac{c}{t}\left(\frac{m}{\log m}+\frac{m}{2\log^2 m}-\frac{2c_1\sqrt{tx}}{\log m}\right) \geq \frac{c}{t}\left(\frac{m}{\log m}+\frac{m}{4\log^2 m}\right) \\ & \geq & \frac{c}{t}\left(\frac{m}{\log m}+\frac{m}{4\log^2 m}\right) \geq \frac{cm}{t\log m}+\frac{x}{m}=d,
\end{eqnarray*}
where we used $c = \frac{1}{576c_1^2}$.
$\hfill\square$

\section{Simple Topological Graphs}

In this section, we prove Theorem \ref{simple}. The following statement will be crucial for our purposes.

\begin{theorem}
\label{keylemma}
Let $G = (V,E)$ be a $k$-quasi-planar simple topological graph with $n$ vertices. Suppose that $G$ has an edge that crosses every other edge. Then we have $|E| \leq n\cdot 2^{\alpha(n)^{c'_k}}$, where $\alpha(n)$ denotes the inverse Ackermann function and $c'_k$ is a constant that depends only on $k$.
\end{theorem}

\noindent \textbf{Proof of Theorem \ref{keylemma}.}  Let $k \geq 5$ and $c'_k=  40\cdot 2^{k^2 + 2k}$.  To simplify the presentation, we do not make any attempt to optimize the value of $c'_k$.  Label the vertices of $G$ from 1 to $n$, i.e., let $V=\{1,2,\ldots,n\}$.  Let $e = uv$ be the edge that crosses every other edge in $G$.  Note that $d(u) = d(v) = 1$.

\begin{figure}[h]
\begin{center}
\includegraphics[width=200pt]{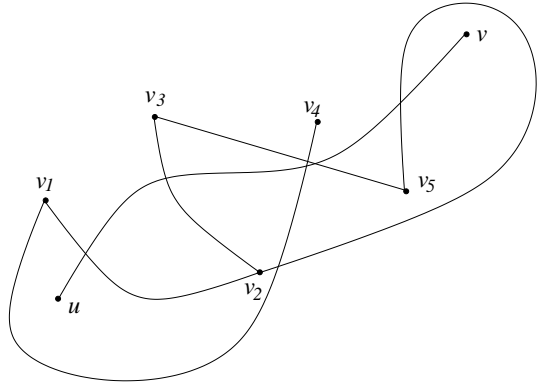}
  \caption{In this example, $S_1 = v_1,v_3,v_4,v_3,v_2$ and $S_2 = v_2,v_2,v_1,v_5,v_5$.}
  \label{s1s2}
 \end{center}
\end{figure}

Let $E'$ denote the set of edges that cross $e$. Suppose without loss of generality that no two of elements of $E'$ cross $e$ at the same point. Let $e_1,e_2,...,e_{|E'|}$ denote the edges in $E'$ listed in the order of their intersection points with $e$ from $u$ to $v$.  We create two sequences of vertices $S_1 = p_1,p_2,...,p_{|E'|}$ and $S_2 = q_1,q_2,...,q_{|E'|}\subset V$, as follows.  For each $e_i \in E'$, as we move along edge $e$ from $u$ to $v$ and arrive at the intersection point with $e_i$, we turn left and move along edge $e_i$ until we reach its endpoint $u_i$.  Then we set $p_i = u_i$.  Likewise, as we move along edge $e$ from $u$ to $v$ and arrive at edge $e_i$, we turn right and move along edge $e_i$ until we reach its other endpoint $w_i$.  Then we set $q_i = w_i$.  Thus, $S_1$ and $S_2$ are sequences of length $|E'|$ over the alphabet $\{1,2,...,n\}$.  See Figure \ref{s1s2} for a small example.

We need two lemmas.  The first one is due to Valtr \cite{valtr}.

\begin{lemma}[Valtr]
\label{regular}
For $l \geq 1$, at least one of the sequences $S_1,S_2$ defined above contains an $l$-regular subsequence of length at least $|E'|/(4l)$. $\hfill\square$
\end{lemma}

Since each edge in $E'$ crosses $e$ exactly once, the proof of Lemma \ref{regular} can be copied almost {\em verbatim} from the proof of Lemma 4 in \cite{valtr} and is left to the reader as an exercise.

For the rest of this section, we set $l = 2^{k^2 + k}$ and $t= 2^k$.

\begin{lemma}
\label{key}
Neither of the sequences $S_1$ and $S_2$ has a subsequence of type $up(l, t)$.
\end{lemma}

\noindent \textbf{Proof.}  By symmetry, it suffices to show that $S_1$ does not contain a subsequence of type $up(l,t)$.  The argument is by contradiction.  We will prove by induction on $k$ that the existence of such a sequence would imply that $G$ has $k$ pairwise crossing edges.  The base cases $k = 1,2$ are trivial.  Now assume the statement holds up to $k-1$.  Let

$$S = s_1,s_2,...,s_{lt}$$

\noindent be our $up(l,t)$ sequence of length $lt$ such that the first $l$ terms are pairwise distinct and for $i = 1,2,...,l$ we have

$$s_i = s_{i + l }  = s_{i + 2l } = s_{i + 3l}= \cdots = s_{i + (t - 1)l}.$$

\noindent For each $i = 1,2,...,l$, let $v_i\in V$ denote the vertex $s_i$.  Moreover, let $a_{i,j}$ be the arc emanating from vertex $v_i$ to the edge $e$ corresponding to $s_{i + jl}$ for $j = 0,1,2,...,t - 1$.  We will think of $s_{i + jl}$ as a point on $a_{i,j}$ very close but not on edge $e$.  For simplicity, we will let $s_{lt + q} = s_q$ for all $ q\in \mathbb{N}$ and  $a_{i,j} = a_{i,j'}$ for all $j  \in \mathbb{Z}$, where $j' \in \{0,1,2,\ldots,t-1\}$ is such that $j \equiv j'$ (mod $t$).  Hence there are $l$ distinct vertices $v_1,...,v_{l}$, each vertex of which has $t$ arcs emanating from it to the edge $e$.

\begin{figure}[h]
\begin{center}
\includegraphics[width=200pt]{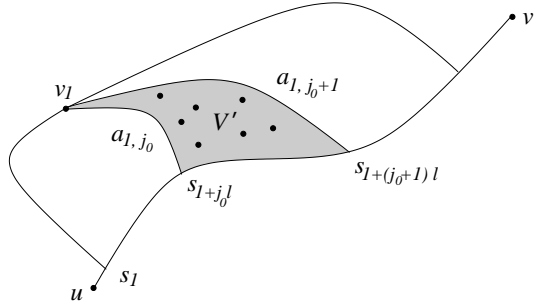}
  \caption{Vertices of $V'$ lie in the region enclosed by $a_{1,j_0},a_{1,j_0 + 1},e$.}
  \label{between}
 \end{center}
\end{figure}

Consider the arrangement formed by the $t$ arcs emanating from $v_1$ and the edge $e$.  Since $G$ is simple, these arcs partition the plane into $t$ regions. By the pigeonhole principle, there is a subset $V' \subset \{v_1,...,v_{l}\}$ of size

$$\frac{l - 1}{t} = \frac{2^{k^2 + k} - 1}{2^k}$$

\noindent such that all of the vertices of $V'$ lie in the same region.  Let $j_0 \in\{0,1,2,...,t-1\}$ be an integer such that $V'$ lies in the region bounded by $a_{1,j_0},a_{1,j_0 + 1},$ and $e$.  See Figure \ref{between}.  In the case $j_0 = t - 1$, the set $V'$ lies in the unbounded region.

Let $v_i \in V'$ and $a_{i,j_0 + j_1}$ be an arc emanating from $v_i$ for $j_1 \geq 1$.  Notice that $a_{i,j_0 + j_1}$ cannot cross both $a_{1,j_0}$ and $a_{1,j_0+1}$, since $G$ is a {\em simple} topological graph.  Suppose that $a_{i,j_0 + j_1}$ crosses $a_{1,j_0 + 1}$.  Then all arcs emanating from $v_i$,

 $$A = \{a_{i,j_0 + 1}, a_{i,j_0 + 2},...,a_{i,j_0 + j_1-1}\}$$

\noindent must also cross $a_{1,j_0 + 1}$.  Indeed, let $\gamma$ be the simple closed curve created by the arrangement

 $$a_{i,j_0 + j_1}\cup a_{1,j_0 + 1}\cup e.$$

\noindent Since $a_{i,j_0 + j_1},a_{1,j_0 + 1},$ and $e$ pairwise intersect at precisely one point, $\gamma$ is well defined.  We define points $x = a_{i,j_0 + j_1} \cap a_{1,j_0  + 1}$ and $y = a_{1,j_0 + 1}\cap e$, and orient $\gamma$ in the direction from $x$ to $y$ along $\gamma$.

\begin{figure}[h]
    \label{gamma}
  \centering
  \subfigure[The case when $ j_0 + j_1 \mod t \leq t - 1$. ]{\label{j1p}\includegraphics[width=0.47\textwidth]{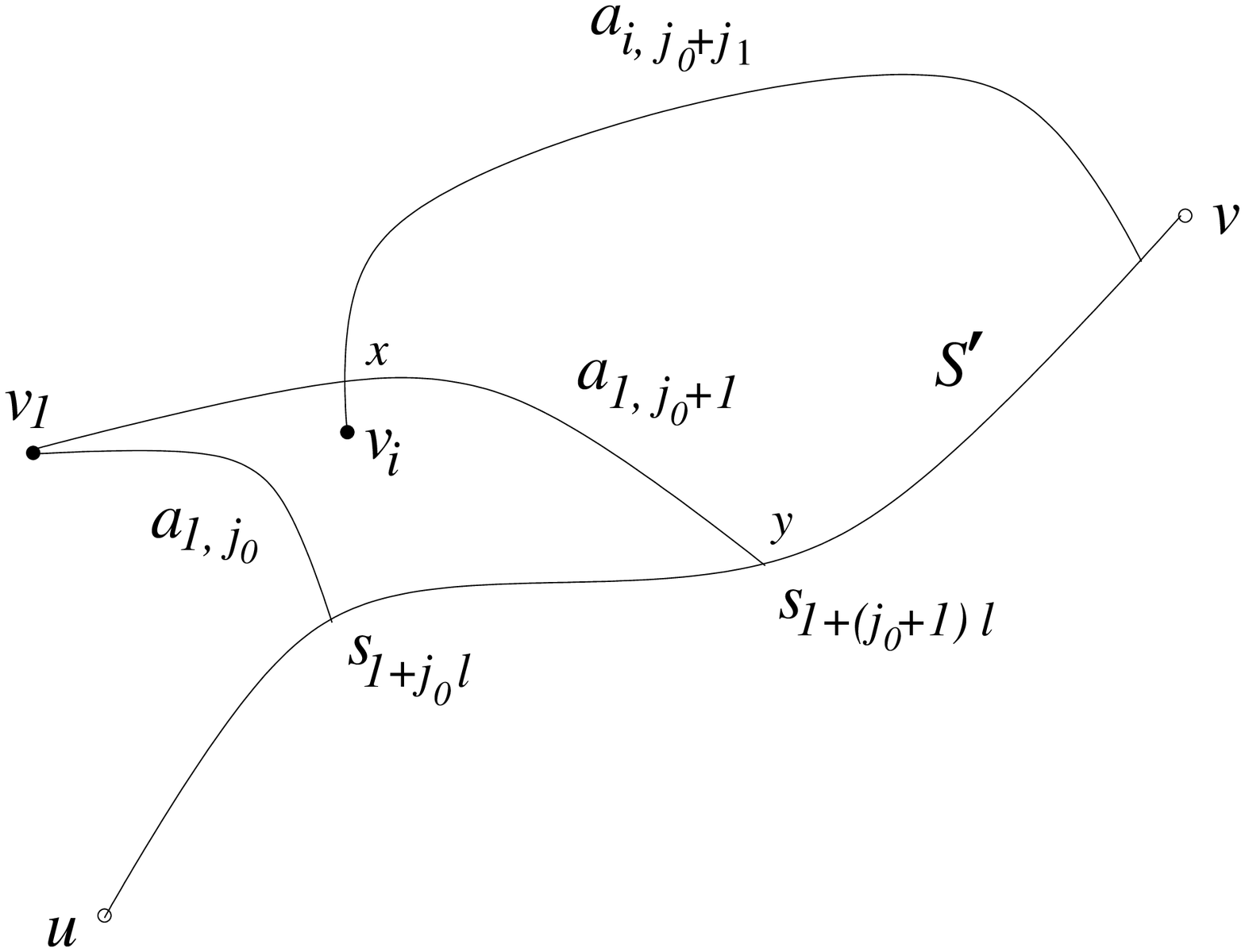}} \hspace{.5cm}
\subfigure[$\gamma$ defined from Figure \ref{j1p}.]{\label{j1}\includegraphics[width=0.47\textwidth]{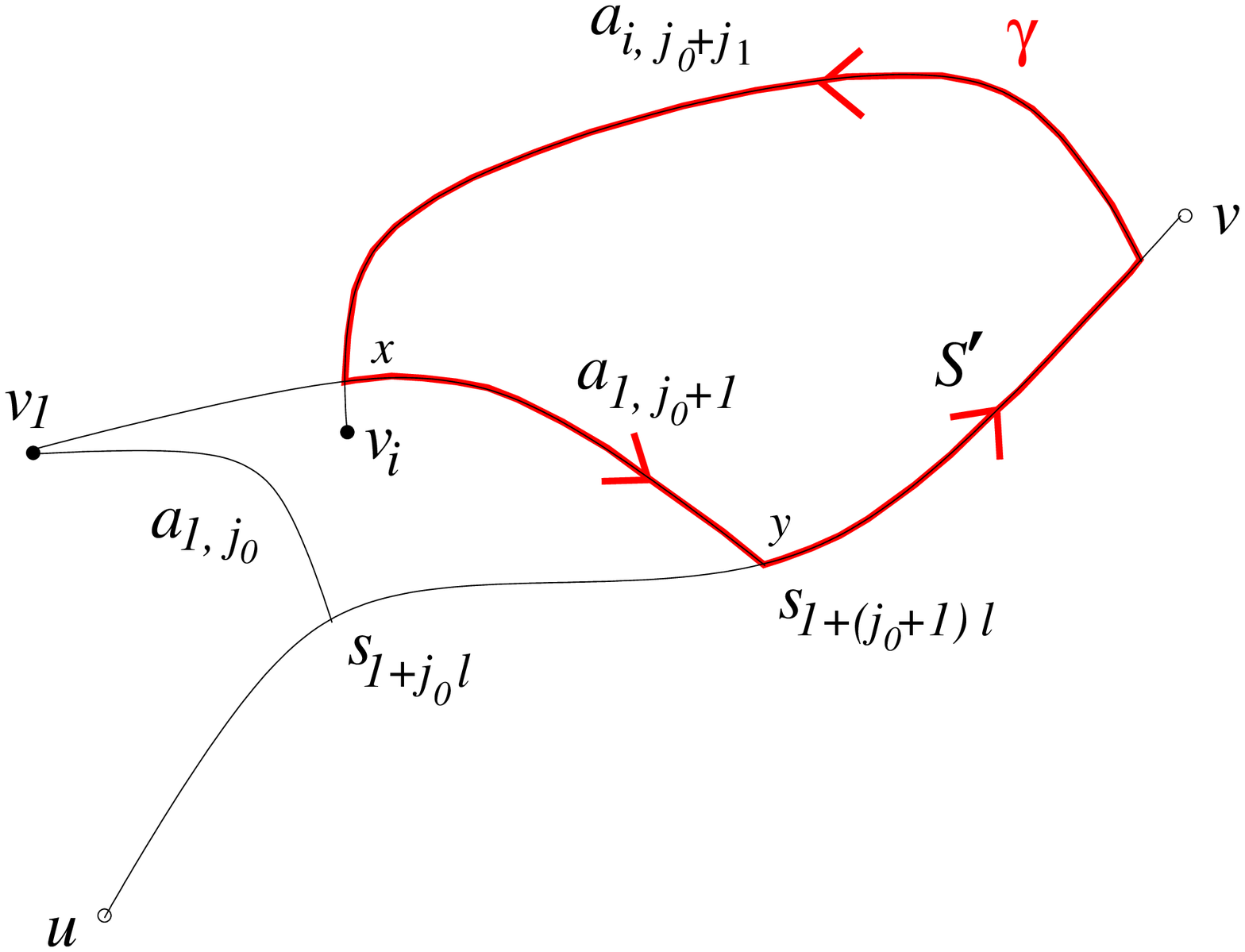}}
  \subfigure[The case when $ j_0 + j_1 \mod t  < j_0$.  Recall $a_{i,j_0 + j_1} = a_{i,j_0+ j_1\mod 2^k}$.]{\label{j13p}\includegraphics[width=0.49\textwidth]{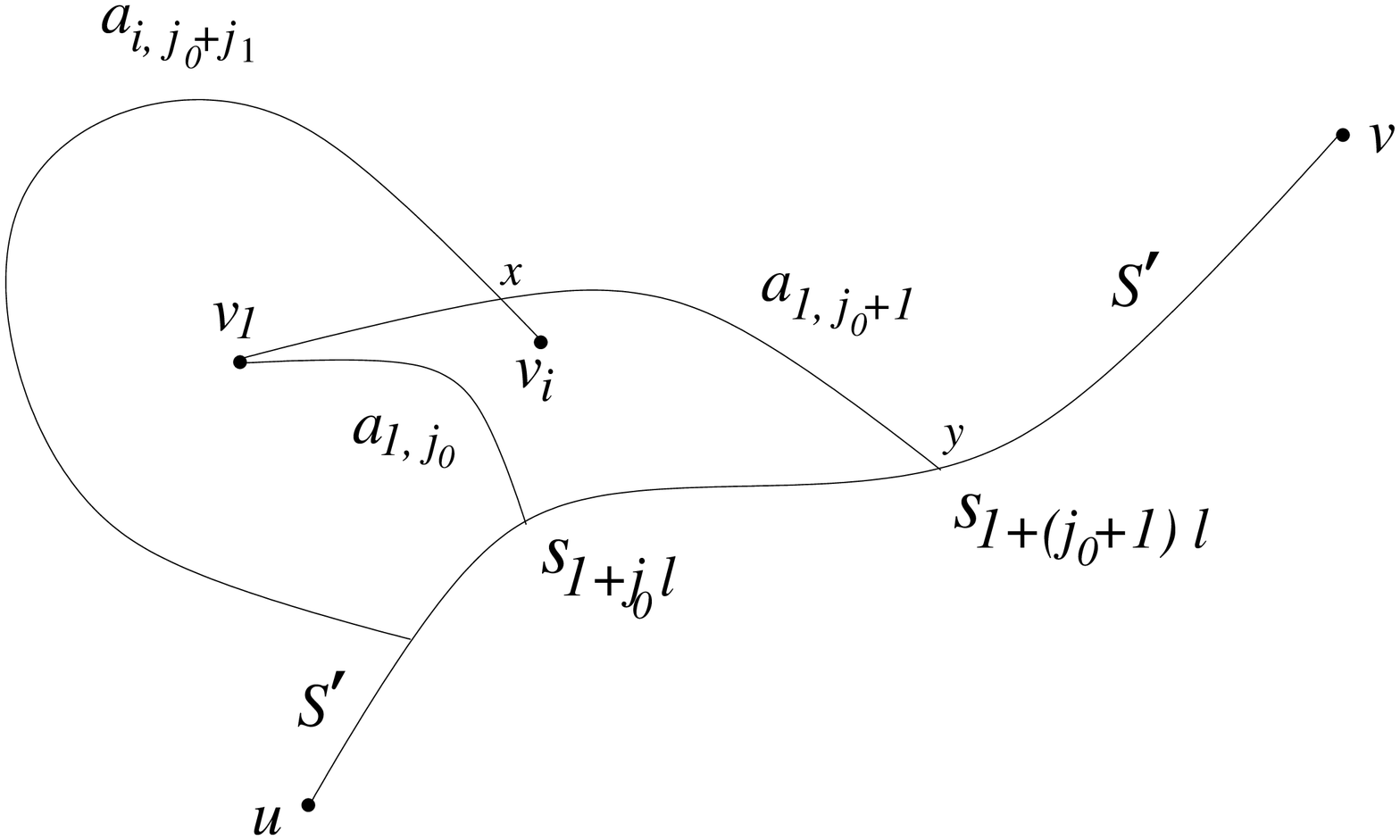}}
\subfigure[$\gamma$ defined from Figure \ref{j13p}.]{\label{j13}\includegraphics[width=0.49\textwidth]{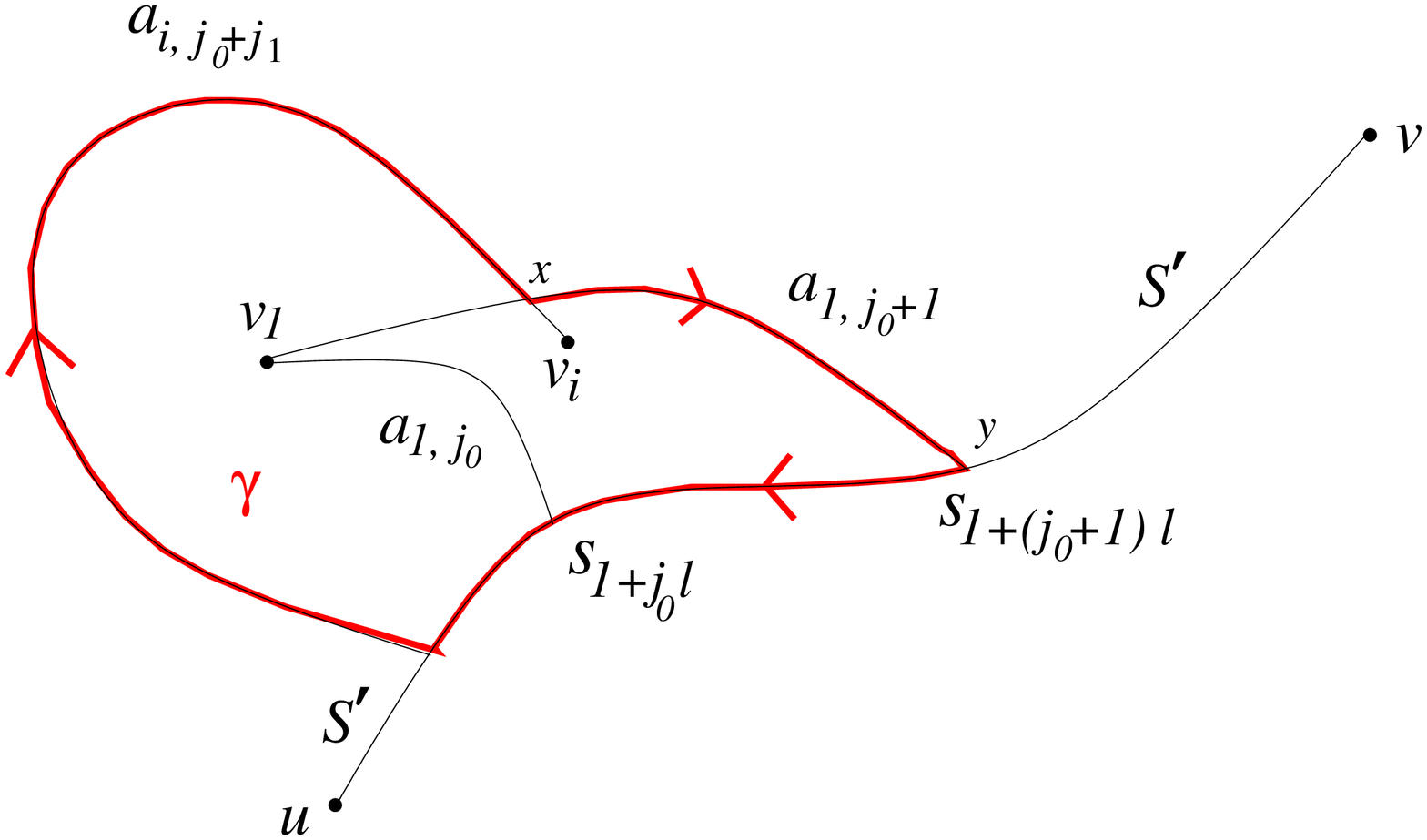}}
                        \caption{Defining $\gamma$ and its orientation.}
\end{figure}

In view of the fact that $a_{i,j_0 + j_1}$ intersects $a_{1,j_0 + 1}$, the vertex $v_i$ must lie to the right of $\gamma$.  Moreover, since the arc from $x$ to $y$ along $a_{1,j_0 + 1}$ is a subset of $\gamma$, the points corresponding to the subsequence

$$ S' = \{s_q \in S \hspace{.2cm}|\hspace{.2cm} 2 + (j_0 + 1)l \leq q \leq (i-1) + (j_0 + j_1)l\}$$

\noindent must lie to the left of $\gamma$.  Hence, $\gamma$ separates vertex $v_i$ and the points of $S'$.  Therefore, using again that $G$ is simple, each arc from $A$ must cross $a_{1,j_0 + 1}$ (these arcs cannot cross $a_{i,j_0 + j_1}$).  See Figure \ref{gamma}.

By the same argument, if the arc $a_{i,j_0 - j_1}$ crosses $a_{1,j_0}$ for $j_1 \geq 1$, then the arcs emanating from $v_i$,

$$a_{i,j_0 - 1},a_{i,j_0 - 2},...,a_{i,j_0 - j_1 + 1}$$

\noindent must also cross $a_{1,j_0}$.  Since $a_{i,j_0 + t/2} = a_{i,j_0 - t/2}$, we have the following observation.

\begin{observation}
\label{obs}
For half of the vertices $v_i \in V'$, the arcs emanating from $v_i$ satisfy
\begin{enumerate}

\item $a_{i,j_0 + 1}, a_{i,j_0 + 2},...,a_{i,j_0 + t/2}$ all cross $a_{1,j_0 + 1}$, or

\item $a_{i,j_0 - 1}, a_{i,j_0 - 2},...,a_{i,j_0 - t/2}$ all cross $a_{1,j_0}$.

\end{enumerate}
\end{observation}

\noindent Since $t/2 = 2^{k-1}$ and

$$\frac{|V'|}{2} \geq \frac{l - 1}{2t} =  \frac{2^{k^2 + k} - 1}{2\cdot 2^{k}} \geq 2^{(k-1)^2 + (k-1)},$$

\noindent by Observation \ref{obs}, we obtain a $up(2^{(k-1)^2 + (k-1)},2^{k - 1})$ sequence such that the corresponding arcs all cross either $a_{1,j_0}$ or $a_{1,j_0 +1}$.  By the induction hypothesis, it follows that there exist $k$ pairwise crossing edges. $\hfill\square$

\medskip

Now we are ready to complete the proof of Lemma \ref{keylemma}.  By Lemma \ref{regular} we know that, say, $S_1$ contains an $l$-regular subsequence of length $|E'|/(4l)$.  By Theorem \ref{klazar} and Lemma \ref{key}, this subsequence has length at most

$$n\cdot l 2^{(lt-3)} \cdot (10l)^{10\alpha(n)^{lt}}.$$

\noindent Therefore, we have

$$ \frac{|E'|}{4\cdot l} \leq n\cdot l 2^{(lt-3)} \cdot (10l)^{10\alpha(n)^{lt}},$$

\noindent which implies

$$|E'| \leq 4 \cdot n\cdot l^2 2^{(lt-3)} \cdot (10l)^{10\alpha(n)^{lt}}.$$

\noindent Since $c'_k = 40\cdot lt = 40\cdot 2^{k^2 + 2k}$, $\alpha(n) \geq 2$ and $k \geq5$, we have

$$|E|=|E'|+1 \leq n\cdot 2^{\alpha^{c'_k}(n)},$$

\noindent which completes the proof of Lemma \ref{keylemma}. $\hfill\square$

\medskip

Now we are in a position to prove Theorem \ref{simple}.

\medskip

\noindent \textbf{Proof of Theorem \ref{simple}.}  Let $G = (V,E)$ be a $k$-quasi-planar simple topological graph on $n$ vertices.  By Lemma \ref{decompose}, there is a subset $E' \subset E$ such that $|E'| \geq c|E|/\log|E|,$ where $c$ is an absolute constant and $E'$ is decomposable.  Hence, there is a partition

$$E' = E_1\cup E_2\cup\cdots \cup E_t$$

\noindent such that each $E_i$ has an edge $e_i$ that intersects every other edge in $E_i$, and for $i\neq j$, the edges in $E_i$ are disjoint from the edges in $E_j$.  Let $V_i$ denote the set of vertices that are the endpoints of the edges in $E_i$, and let $n_i = |V_i|$.  By Lemma \ref{keylemma}, we have

$$|E_i| \leq n_i 2^{\alpha^{c'_k}(n_i)} + 2n_i,$$

\noindent where the $2n_i$ term accounts for the edges that share a vertex with $e_i$.  Hence,

$$\frac{c|E|}{\log |E|} \leq \sum\limits_{i = 1}^t n_i 2^{\alpha^{c'_k}(n_i)} + 2n_i \leq n2^{\alpha^{c'_k}(n)} + 2n,$$

\noindent Therefore, we obtain

$$|E| \leq (n\log n) 2^{\alpha^{c_k}(n)},$$

\noindent for a sufficiently large constant  $c_k$.
$\hfill\square$

\section{$x$-Monotone Topological Graphs}

The aim of this section is to prove Theorem \ref{xmono}.

\medskip

\noindent \textbf{Proof of Theorem \ref{xmono}.}   For $k \geq 2$, let $g_k(n)$ be the maximum number of edges in a $k$-quasi-planar topological graph whose edges are drawn as $x$-monotone curves.  We will prove by induction on $n$ that

  $$g_k(n) \leq 2^{c k^6}n\log n$$

\noindent where $c$ is a sufficiently large absolute constant.

The base case is trivial.  For the inductive step, let $G = (V,E)$ be a $k$-quasi-planar topological graph whose edges are drawn as $x$-monotone curves, and let the vertices be labeled $1,2,...,n$.  Let $L$ be a vertical line that partitions the vertices into two parts, $V_1$ and $V_2$, such that $|V_1| = \lfloor n/2\rfloor$ vertices lie to the left of $L$, and $|V_2| = \lceil n/2\rceil$ vertices lie to the right of $L$.  Furthermore, let $E_1$ denote the set of edges induced by $V_1$, let $E_2$ denote the set of edges induced by $V_2$, and let $E'$ be the set of edges that intersect $L$. Clearly, we have

$$|E_1| \leq g_k(\lfloor n/2\rfloor) \hspace{0.5cm}\textnormal{and}\hspace{0.5cm} |E_2| \leq g_k(\lceil n/2\rceil).$$

\noindent It suffices that show that

\begin{equation}
\label{goal}
|E'|\leq 2^{c k^6/2}n,
\end{equation}

\noindent since this would imply

$$g_k(n) \leq g_k(\lfloor n/2\rfloor) + g_k(\lceil n/2\rceil)+ 2^{c k^6/2}n \leq  2^{c k^{6}} n\log n. $$

In the rest of the proof, we only consider the edges belonging to $E'$.  For each vertex $v_i \in V_1$, consider the graph $G_i$ whose vertices are the edges with $v_i$ as a left endpoint, and two vertices in $G_i$ are adjacent if the corresponding edges cross at some point to the left of $L$.  Since $G_i$ is an \emph{incomparability graph} (see \cite{dilworth}, \cite{fox2}) and does not contain a clique of size $k$, $G_i$ contains an independent set of size $|E(G_i)|/(k-1)$.  We keep all edges that correspond to the elements of this independent set, and discard all other edges incident to $v_i$.  After repeating this process for all vertices in $V_1$, we are left with at least $|E'|/(k-1)$ edges.

Now we continue this process on the other side.  For each vertex $v_j \in V_2$, consider the graph $G_j$ whose vertices are the edges with $v_j$ as a right endpoint, and two vertices in $G_j$ are adjacent if the corresponding edges cross at some point to the right of $L$.  Since $G_j$ is an incomparability graph and does not contain a clique of size $k$, $G_j$ contains an independent set of size $|E(G_j)|/(k-1)$.  We keep all edges that corresponds to this independent set, and discard all other edges incident to $v_j$.  After repeating this process for all vertices in $V_2$, we are left with at least $|E'|/(k-1)^2$ edges.

We order the remaining edges $e_1,e_2,...,e_m$ in the order in which they intersect $L$ from bottom to top. (We assume without loss of generality that any two intersection points are distinct.)  Define two sequences, $S_1 = p_1,p_2,...,p_m$ and $S_2 = q_1,q_2,...,q_m$, such that $p_i$ denotes the left endpoint of edge $e_i$ and $q_i$ denotes the right endpoint of $e_i$.  We need the following lemma.

\begin{lemma}
\label{key2}
Neither of the sequences $S_1$ and $S_2$ has subsequence of type up-down-up$(k^3+2)$.
\end{lemma}

\noindent \textbf{Proof.}  By symmetry, it suffices to show that $S_1$ does not have a subsequence of type \emph{up-down-up}$(k^3 + 2)$.  Suppose for contradiction that $S_1$ does contain such a subsequence.  Then there is a sequence

$$S = s_1,s_2,...,s_{3(k^3 + 2)-2}$$

\noindent such that the integers $s_1,...,s_{k^3 + 2}$ are pairwise distinct and

$$s_i = s_{2(k^3 + 2) - i} = s_{2(k^3 + 2) - 2 + i}$$

for $i = 1,2,...,k^3 + 2$.

For each $i\in \{1,2,...,k^3 + 2\}$, let $v_i \in V_1$ denote the label (vertex) of $s_i$ and let $x_i$ denote the $x$-coordinate of the vertex $v_i$.  Moreover, let $a_i$ be the arc emanating from vertex $v_i$ to the point on $L$ that corresponds to $s_{2(k^3+2) - i}$. Let $A = \{a_2,a_3,...,a_{k^3 + 1}\}$. Note that the arcs in $A$ are enumerated downwards with respect to their intersection points with $L$, and they correspond to the elements of the ``middle" section of the up-down-up sequence.  We define two partial orders on $A$ as follows.

$$\begin{array}{cccccc}
   a_i \prec_1 a_j & \textnormal{if} &  i < j, &x_i < x_j & \textnormal{and the arcs $a_i,a_j$ do not intersect,}\\\\
   a_i \prec_2 a_j & \textnormal{if} &  i < j, &x_i > x_j & \textnormal{and the arcs $a_i,a_j$  do not intersect.}
 \end{array}$$

Clearly, $\prec_1$ and $\prec_2$ are partial orders.  If two arcs are not comparable by either $\prec_1$ or $\prec_2$, then they must cross.  Since $G$ does not contain $k$ pairwise crossing edges, by Dilworth's theorem, there exist $k$ arcs $\{a_{i_1},a_{i_2},...,a_{i_k}\}$ such that they are pairwise comparable by either $\prec_1$ or $\prec_2$.  Now the proof falls into two cases.

\medskip

\begin{figure}[h]
  \centering
\subfigure{ \includegraphics[width=0.3 \textwidth]{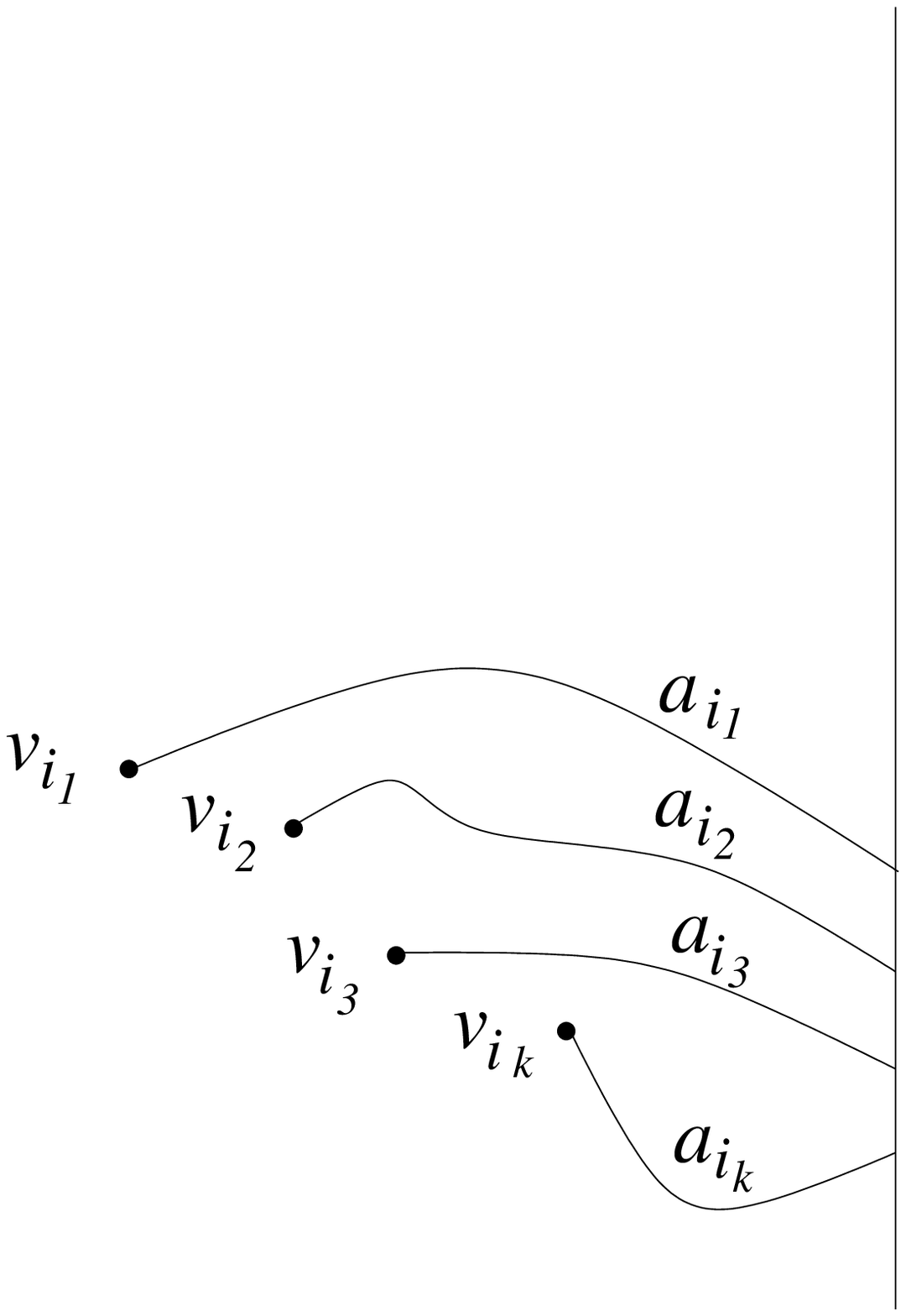}} \hspace{2cm}
\subfigure{ \includegraphics[width=0.38\textwidth]{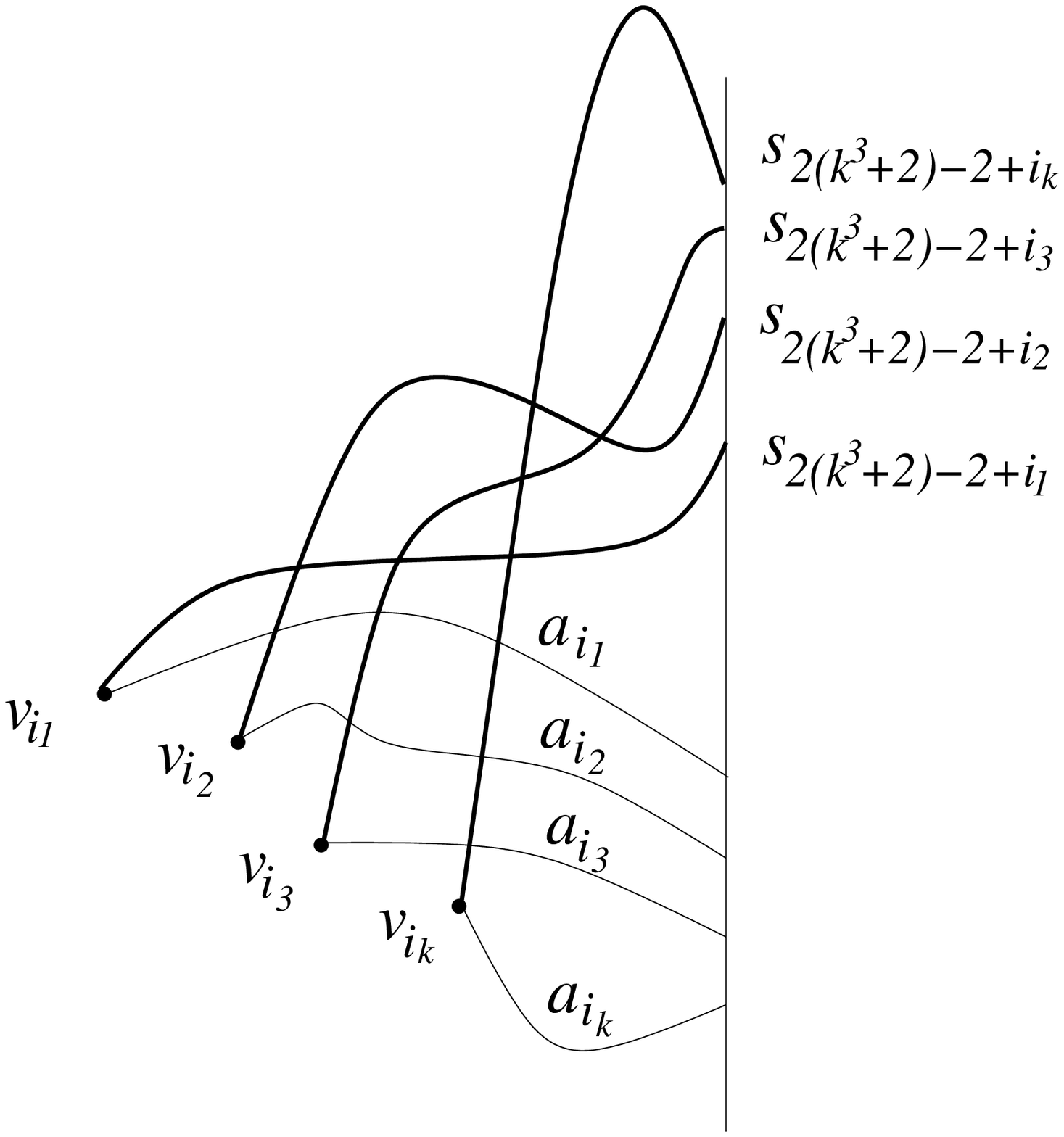}}
                        \caption{Case 1.}
  \label{prec1}
\end{figure}

\noindent \emph{Case 1.}  Suppose that $a_{i_1} \prec_1 a_{i_2} \prec_1 \cdots \prec_1 a_{i_k}$.  Then the arcs emanating from $v_{i_1},v_{i_2},...,v_{i_k}$ to the points corresponding to $s_{2(k^3+2)-2 + i_1},s_{2(k^3 + 2)-2 + i_2},...,s_{2(k^3 + 2)-2 + i_k}$ are pairwise crossing.  See Figure \ref{prec1}.

\medskip

\begin{figure}[h]
  \centering
\subfigure{ \includegraphics[width=0.3\textwidth]{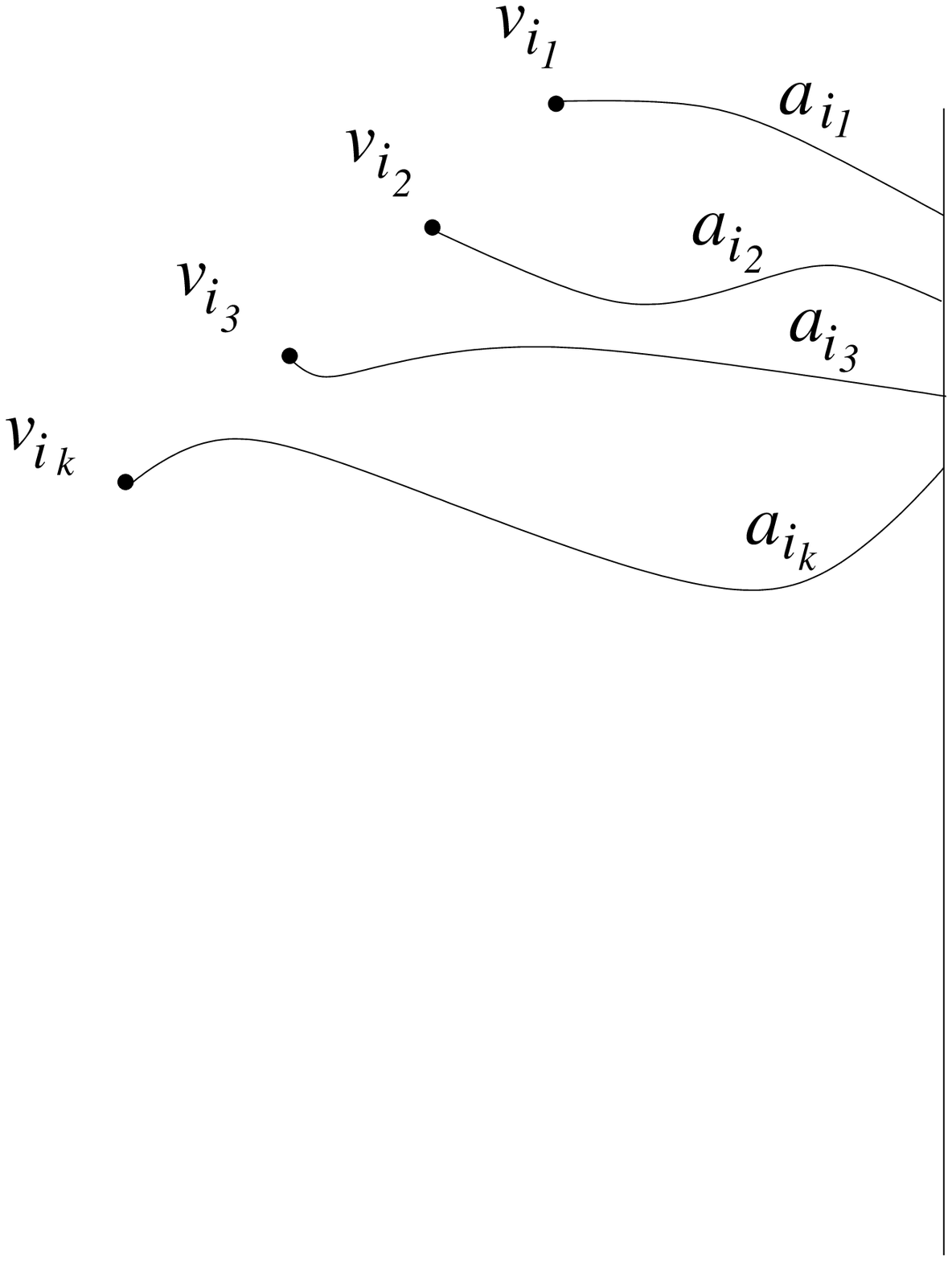}} \hspace{2cm}
\subfigure{ \includegraphics[width=0.38\textwidth]{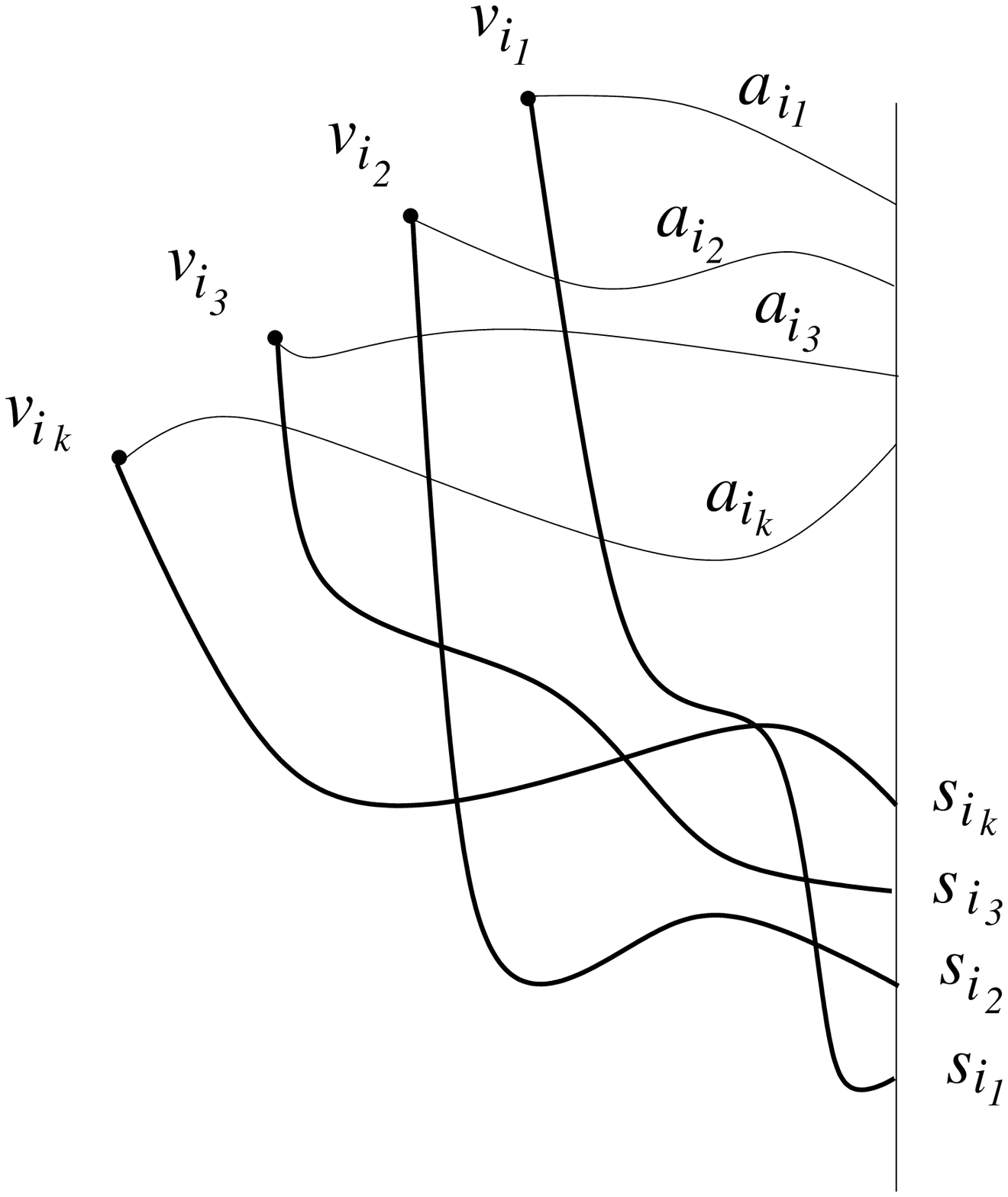}}
                        \caption{Case 2.}
  \label{prec2}
\end{figure}

\noindent \emph{Case 2.}    Suppose that $a_{i_1} \prec_2 a_{i_2} \prec_2 \cdots \prec_2 a_{i_k}$. Then the arcs emanating from $v_{i_1},v_{i_2},...,v_{i_k}$ to the points corresponding to $s_{i_1},s_{i_2},...,s_{i_k}$ are pairwise crossing.  See Figure \ref{prec2}.

$\hfill\square$

\medskip

We are now ready to complete the proof of Theorem \ref{xmono}.  By Lemma \ref{regular}, we know that, $S_1$, say, contains a $(k^3+2)$-regular subsequence of length

$$\frac{|E'|}{4(k^3 + 2)(k-1)^2}.$$

\noindent By Lemmas \ref{updownup} and \ref{key2}, this subsequence has length at most $2^{c' k^6}n$, where $c'$ is an absolute constant.  Hence, we have

$$\frac{|E'|}{4(k^3 + 2)(k-1)^2} \leq 2^{c' k^6}n,$$

\noindent which implies that

$$|E'| \leq 4k^52^{c'k^6}n \leq 2^{c k^6/2} n$$

\noindent for a sufficiently large absolute constant $c$.
$\hfill\square$

\newpage

\bibliographystyle{plain}

\end{document}